\documentclass[]{amsart}
\usepackage[foot]{amsaddr}
\usepackage[a4paper,top=2.5cm,bottom=3cm,inner=3cm,outer=3cm]{geometry}

\usepackage[utf8]{inputenc}

\usepackage{amsmath,amsthm}
\usepackage{amsfonts,amssymb,mathrsfs}

\usepackage{xy}
\usepackage[bookmarks=false,breaklinks=true]{hyperref}
\usepackage{float}
\usepackage{graphicx}
\usepackage[percent]{overpic}
\usepackage{color}

\usepackage{epigraph}
\definecolor{myurlcolor}{rgb}{0,0,0.7}
\usepackage{hyperref}
\hypersetup{colorlinks,
linkcolor=myurlcolor,
citecolor=myurlcolor,
urlcolor=myurlcolor}
\usepackage[mathscr]{euscript}

\input xy
\xyoption{all}

\newcommand{\maps}{\colon}    

\newcommand{\Z}{{\mathbb Z}}  
\newcommand{\R}{{\mathbb R}}  

\newcommand{\SO}{{\rm SO}}      
\newcommand{\D}{\mathrm{D}}     
\newcommand{\E}{\mathrm{E}}     

\newcommand{\define}[1]{{\bf \boldmath{#1}}}

\theoremstyle{definition}

        \newcommand{\be}{\begin{equation}}
        \newcommand{\ee}{\end{equation}}
        \newcommand{\ba}{\begin{eqnarray}}
        \newcommand{\ea}{\end{eqnarray}}
        \newcommand{\ban}{\begin{eqnarray*}}
        \newcommand{\ean}{\end{eqnarray*}}
        \newcommand{\barr}{\begin{array}}
        \newcommand{\earr}{\end{array}}

\begin{document}
\title{The Icosidodecahedron }
\author[Baez]{John C.\ Baez} 
\address{Department of Mathematics, University of California, Riverside CA, 92521, USA}
\date{\today}
\maketitle

\begin{center}
\includegraphics[width = 15em]{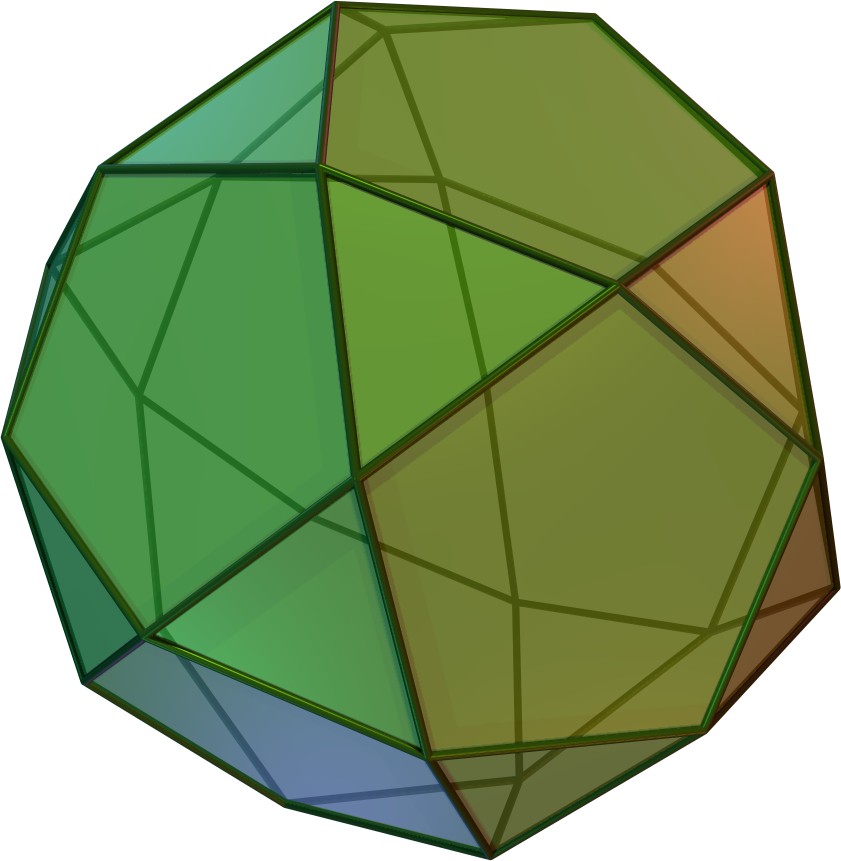} 
\end{center}

\vskip 1em
The icosidodecahedron can be built by truncating either a regular icosahedron or a regular dodecahedron. It has 30 vertices, one at the center of each edge of the icosahedron---or equivalently, one at the center of each edge of a dodecahedron.  It is a beautiful, highly symmetrical shape.  

But the icosidodecahedron is just a shadow of a more symmetrical shape with twice as many vertices, which lives in a space with twice as many dimensions!  Namely, it is a projection down to 3d space of a 6-dimensional polytope with 60 vertices.  Even better,  it is also a slice of a still more symmetrical 4d polytope with 120 vertices, which in turn is the projection down to 4d space of a even more symmetrical 8-dimensional polytope with 240 vertices.  Note how the numbers keep doubling: 30, 60, 120 and 240.

To understand all this, start with the group of rotational symmetries of the icosahedron.  This is a 60-element subgroup of the rotation group $\SO(3)$.  The group of unit quaternions is a double cover of $\SO(3)$, so this 60-element subgroup has a double cover, called the \define{binary icosahedral group}, consisting of 120 unit quaternions.   With a suitable choice of coordinates, we can take these to be
\[ \pm 1 , \quad \frac{\pm 1 \pm i \pm j \pm k}{2}, \quad \frac{\pm i \pm \phi j \pm \Phi k}{2} \]
together with everything obtained from these by even permutations of $1, i, j,$ and $k$, where
\[  \phi = \frac{\sqrt{5} - 1}{2}, \quad \Phi = \frac{\sqrt{5} + 1}{2}, \]
are the `little' and `big' golden ratios, respectively.  These 120 unit quaternions are the vertices of a convex polytope in 4 dimensions.  In fact this is a regular polytope, called the \define{600-cell} since it has 600 regular tetrahedra as faces \cite{Coxeter}.  

Here is the 600-cell projected down to 3d space, as drawn using Robert Webb's Stella software \cite{Webb}:

\begin{center}
\includegraphics[width = 18em]{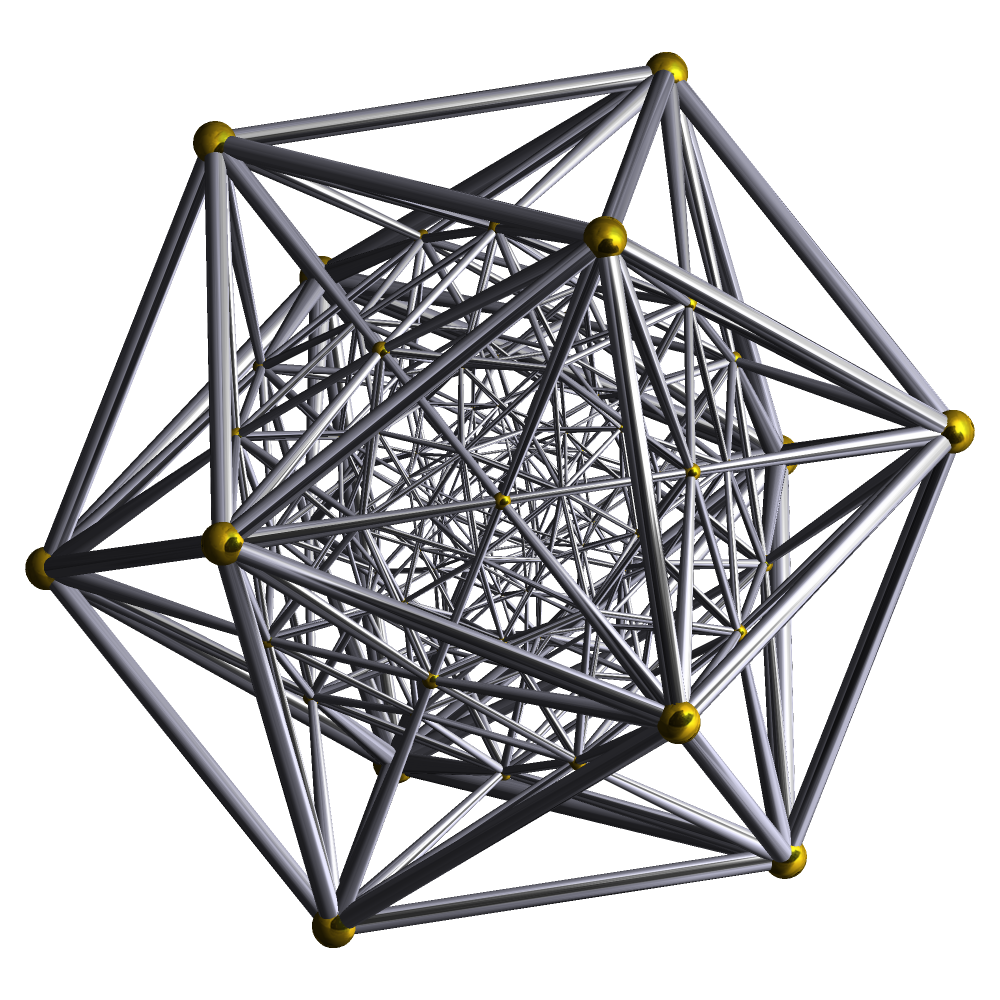} 
\end{center}

If we slice the 600-cell halfway between two of its opposite vertices, we get an icosidodecahedron.   This is easiest to see by intersecting the 600-cell with the space of \define{pure imaginary quaternions}:
\[  \{ ai + bj + ck : \; a,b,c \in \mathbb{R} \}  \]
Of the 600-cell's vertices, those that lie in this 3-dimensional space are
\[  \pm i, \pm j, \pm k, \]
which are the vertices of an octahedron, and the points
\[   \displaystyle{ \frac{\pm i \pm \phi j \pm \Phi k}{2} ,  \quad
\frac{\pm j \pm \phi k \pm \Phi i}{2} , \quad
\frac{\pm k \pm \phi i \pm \Phi j}{2} 
} 
\]
which are the vertices of three `golden boxes'.  A \define{golden box} is the 3d analogue of a golden rectangle: its three sides are in the proportions $\phi, 1$ and $\Phi.$   It is well-known that the vertices of an octahedron and the golden boxes, taken together in this way, are the vertices of an icosidodecahedron:

\begin{center}
\includegraphics[width = 15em]{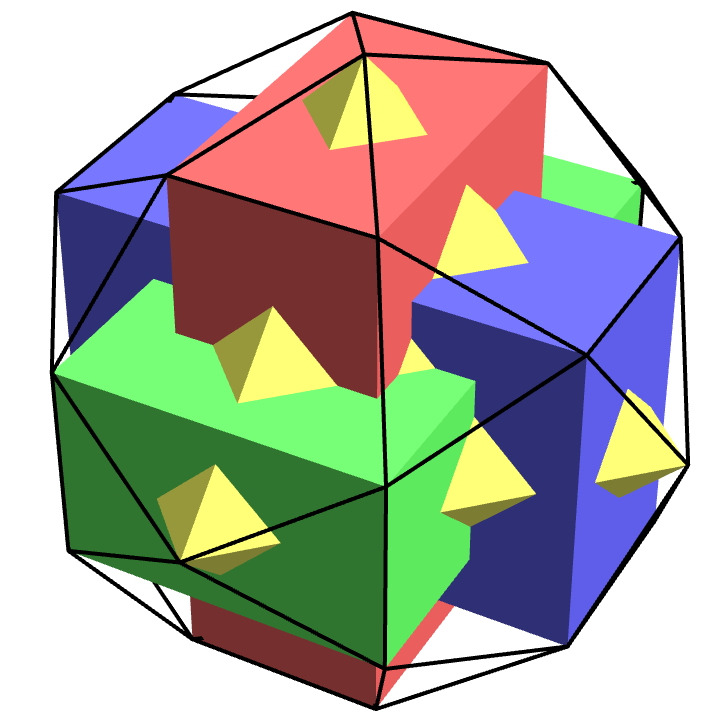} 
\end{center}

But we are not done with the binary icosahedral group---far from it!   Since these 120 quaternions are closed under quaternion multiplication, their integer linear combinations form a subring of the quaternions, which Conway and Sloane \cite{CS} call the \define{icosians}.   Since any icosian can be written as $a + bi + cj + dk$ where the numbers $a,b,c,d \in \R$  are of the form $x + y \sqrt{5}$ with $x,y$ rational, any icosian gives an 8-tuple of rational numbers.  However, we do not get all 8-tuples of rationals this way, only those lying in a certain lattice in $\R^8$.  There is a way to think of this lattice as a rescaled copy of the famous $\E_8$ lattice.  To do this, Conway and Sloane put a new norm on the icosians as follows.  The usual quaternionic norm is
\[    \|a + bi + cj + dk\|^2 = a^2 + b^2 + c^2 + d^2  \]
But for an icosian this norm is always of the form $x + \sqrt{5} y$ for some rationals $x$ and $y$.  Conway and Sloane define a new norm on the icosians by setting
\[        |a + bi + cj + dk|^2 = x + y .\]
With this new norm, Conway and Sloane show the icosians are isomorphic to a rescaled version of a the lattice in $\R^8$ that gives the densest packing of spheres in 8 dimensions: the $\E_8$ lattice \cite{Viazovska}.

The 240 shortest nonzero vectors in the $\E_8$ lattice are the vertices of an 8-dimensional convex polytope called the \define{$\E_8$ root polytope}:

\begin{center}
\includegraphics[width = 18em]{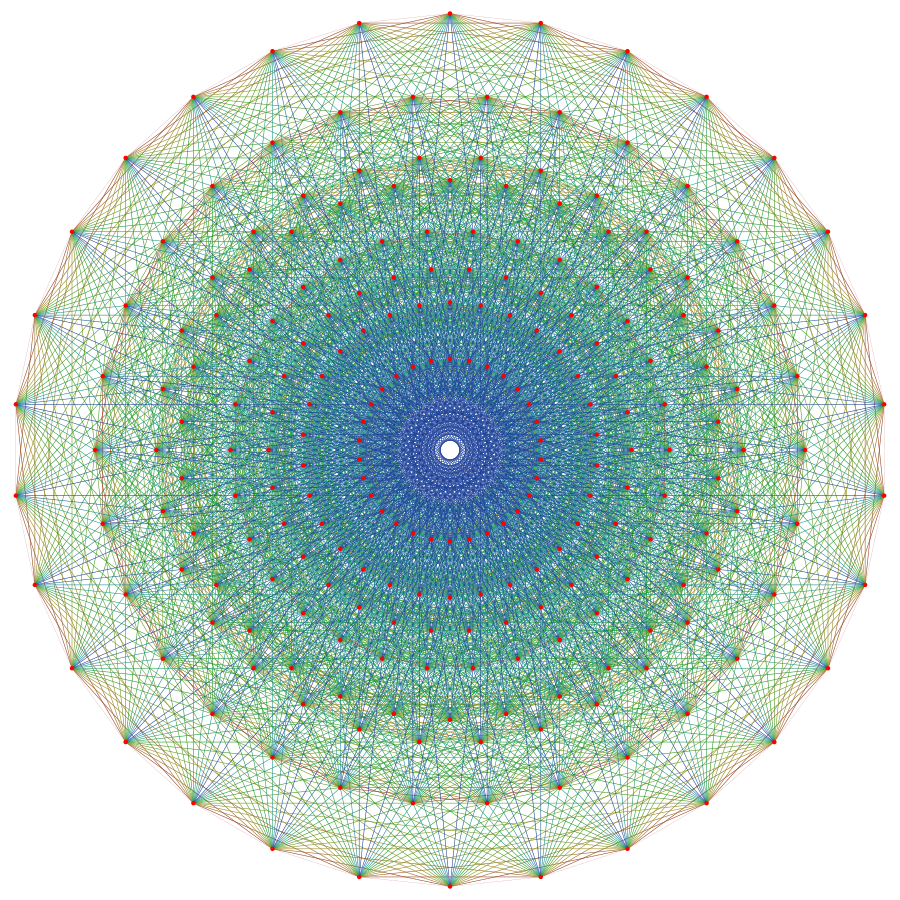} 
\end{center}

\noindent
However, if we remember that each of these 240 vectors came from a quaternion, we can also think of them as 240 quaternions.  These turn out to be the vertices of \emph{two} 600-cells in the quaternions!  In the usual quaternionic norm, one of these 600-cells is larger than the other by a factor of $\Phi$. 

In fact, there is an orthogonal projection from $\R^8$ down to $\R^4$ that maps the $\E_8$ root polytope to the 600-cell.  So, in a very real sense, the 600-cell is the `shadow' of a polytope with twice as many vertices, living in a space whose dimension is twice as large.  And as a spinoff, this fact gives the same sort of relationship between the icosidodecahedron and a 6-dimensional polytope.

The key is to look at the \define{pure imaginary icosians}: icosians of the form $a i + b j + c k $ for real $a,b,c$.  Since $a,b$ and $c$ are each of the form $x + \sqrt{5}y$ with $x$ and $y$ rational, any pure imaginary icosian gives a 6-tuple of rational numbers.  We do not get all 6-tuples of rationals this way, but only those lying in a certain lattice.  We have
\[     \|ai + bj + ck\|^2 = a^2 + b^2 + c^2  \]
For a pure imaginary icosian this is always of the form $x + \sqrt{5} y$ for some rationals $x$ and $y$.  So, we can define a new norm on the pure imaginary icosians by
\[          |ai + bj + ck|^2 = x + y \]
With this new norm, the pure imaginary icosians are isomorphic to a rescaled version of a lattice in $\R^6$ called the \define{$\D_6$ lattice}.

The 60 shortest nonzero vectors in the $\D_6$ lattice are called the \define{roots} of this lattice, and they are the vertices of a 6-dimensional convex polytope called the \define{$\D_6$ root polytope}.   There is an orthogonal projection from $\R^6$ to $\R^3$ that maps this polytope to an icosidodecahedron.   In fact 30 vertices of the $\D_6$ root polytope map to the vertices of this icosidodecahedron, while the other 30 map to vertices of a second, smaller icosidodecahedron.

Let us see some details.  The usual coordinatization of the $\D_6$ lattice in Euclidean $\R^6$ is
\[ \D_6 = \left\{ (x_1, \dots, x_6) : \; x_i  \in \mathbb{Z}, \; \sum_i x_i \in 2\Z \right\} \subset \R^6 .\]
Its roots are the vectors
\[   (\pm 1, \pm 1, 0, 0, 0, 0) \]
and all vectors obtained by permuting the six coordinates.   We shall see that these vectors are sent to the vertices of an icosidodecahedron by the linear map $T \maps  \R^6 \to \R^3$ given as a $3 \times 6$ matrix by
\[  \left( \begin{array}{cccccc}
\Phi &  \Phi  & -1 & -1 & 0 &  0 \\
0 &  0  & \Phi &  -\Phi & -1 & 1 \\
-1 &  1 &  0 &  0 &  \Phi  & \Phi
\end{array} \right)
\]
The rows of this matrix are orthogonal, all with the same norm, so after rescaling it by a constant factor we obtain an orthogonal projection.  The columns of this matrix are six vertices of an icosahedron, chosen so that we never have a vertex and its opposite.   For any pair of columns, they are either neighboring vertices of the icosahedron, or a vertex and the opposite of a neighboring vertex.

The map $T$ thus sends any $\D_6$ root to either the sum or the difference of two neighboring icosahedron vertices.  In this way we obtain all possible sums and differences of neighboring vertices of the icosahedron.    It is easy to see that the sums of neighboring vertices give the vertices of an icosidodecahedron, since by definition the icosidodecahedron has vertices at the midpoints of the edges of a regular icosahedron.  It is less obvious that the differences of neighboring vertices of the icosahedron give the vertices of a second, smaller icosidodecahedron.   But thanks to the symmetry of the situation, we can check this by considering just one example.  In fact the vectors defining the vertices of the larger icosidodecahedron turn out to be precisely $\Phi$ times the vectors defining the vertices of the smaller one!  

The beauties we have just seen are part of an even larger pattern relating all the non-crystallographic Coxeter groups to crystallographic Coxeter groups.  For more, see the work of Fring and Korff \cite{FringKorff}, Boehm, Dechant and Twarock \cite{DBT}, and the many papers they refer to.  Fring and Korff apply these ideas to integrable systems in physics, while the latter authors explore connections to affine Dynkin diagrams.  For more relations between the icosahedron and $\E_8$, see \cite{Baez2}.

\subsubsection*{Acknowledgements}

I thank Greg Egan for help with these ideas.  The image of an icosidodecahedron was created by Cyp and was put on \href{https://commons.wikimedia.org/wiki/File:Icosidodecahedron.jpg}{Wikicommons} with a Creative Commons Attribution-Share Alike 3.0 Unported license.
The 600-cell was made using \href{http://www.software3d.com/Stella.php}{Robert Webb's Stella software} and is available on \href{https://commons.wikimedia.org/wiki/File:Schlegel_wireframe_600-cell_vertex-centered.png}{Wikicommons} with a  Creative Commons Attribution-Share Alike 3.0 Unported license.  The icosidodecahedron with three golden boxes and an octahedron inscribed in it was created by Rahul Narain on \href{https://mathstodon.xyz/web/@narain/109342508956447823}{Mathstodon}.   The projection of the $\E_8$ roots to the plane was created by Claudio Rocchini and put on \href{https://commons.wikimedia.org/wiki/File:E8_graph.svg}{Wikicommons} with a Creative Commons Attribution 3.0 Unported license. 

This article is an expanded of an earlier article \cite{Baez1}.  The only thing I left out here is a rotating image of two icosidodecahedra prepared by Greg Egan and the matrix describing a linear map $S \maps \R^8 \to \R^4$ that when suitably rescaled gives a projection mapping the $\E_8$ lattice in its usual coordinatization
\[ \{ x \in \R^8: \, \textrm{all } x_i \in \Z \textrm{ or all } x_i \in \Z + \frac{1}{2} \textrm{ and } \sum_i x_i \in 2Z \} \]
to the icosians, and thus mapping the 240 $\E_8$ roots to two 600-cells.   For completeness, here is that matrix:
\[ \left( \begin{array}{cccccccc}
\Phi+1 & \Phi -1 & 0  & 0 &  0 &  0 &   0  & 0 \\
0 & 0 & \Phi &  \Phi  & -1 & -1 & 0 & 0 \\
0 & 0  & 0 &  0  & \Phi &  -\Phi & -1 & 1   \\
0 & 0 & -1 &  1 &  0 &  0 &  \Phi  & \Phi
\end{array} \right)
\]

\end{document}